\documentclass[a4paper,11pt]{article}
\usepackage{graphicx}
\usepackage{amssymb}
\usepackage{latexsym, bm}
\usepackage{multicol}
\usepackage{indentfirst}
\usepackage{amssymb,amsfonts}
\usepackage{amsmath}
\usepackage{setspace}
\usepackage{enumerate}
\usepackage{CJK}

\newtheorem{theorem}{Theorem}
\newtheorem{claim}{Claim}
\newtheorem{lemma}{Lemma}
\newtheorem{corollary}{Corollary}

\newtheorem{definition}{Definition}

\newtheorem{fact}{Fact}

\usepackage{amssymb}
\begin{spacing}{1.0}
\begin{document}
\title{Ramsey numbers of large  even cycles and fans}

\date{}
%\footnote{Center for Discrete Mathematics, Fuzhou University,
%Fuzhou, 350108, P.~R.~China.  Email: {\em chunlin\_you@163.com, linqizhong@fzu.edu.cn}  }}
%
%\author{Chunlin You\;\; \text{and} \;\;
%Qizhong Lin\footnote{Corresponding author. Supported in part by NSFC(No.\ 12171088).}
%}
%\maketitle
\author{{Chunlin You\footnote{school of mathematics and statistics, Yancheng Teachers University, Yancheng
224002, P.~R.~China. Email: {\tt chunlin\_you@163.com}.}
~and
 {Qizhong Lin\footnote{Center for Discrete Mathematics, Fuzhou University,
Fuzhou, 350108 P.~R.~China. Email: {\tt linqizhong@fzu.edu.cn}. Supported by the NSFC grant
(No.\ 11671088) }}}}
\maketitle

\begin{abstract}
For graphs $F$ and $H$, the Ramsey number $R(F, H)$ is the smallest positive integer $N$ such that any red/blue edge coloring of $K_N$ contains either a red $F$ or a blue $H$.	 Let $C_n$ be a cycle of length $n$  and $F_n$ be a fan consisting of $n$ triangles all sharing a common vertex.
In this paper, we prove that for all sufficiently large $n$,
\[
R(C_{2\lfloor an\rfloor}, F_n)= \left\{ \begin{array}{ll}
(2+2a+o(1))n & \textrm{if $1/2\leq a< 1$,}\\
(4a+o(1))n & \textrm{if $ a\geq 1$.}
\end{array} \right.
\]

\medskip
		
		
%\medskip
		
%{\bf Mathematics Subject Classification:} \  05C55;\;\;05D10
		
\end{abstract}

\section{Introduction}
For graphs $H_1$ and $H_2$, the Ramsey number $r(H_1,H_2)$ is defined as the smallest integer $N$ such that for any red/blue edge coloring of $K_N$, there exists either a red $H_1$ or a blue $H_2$. The existence of Ramsey number $R(H_1,H_2)$ follows from Ramsey \cite{ram}.

Let $C_n$ and $K_n$ be a cycle and a complete graph on $n$ vertices, respectively.
A fan $F_n$ is a graph on $2n+1$ vertices with a vertex $v$, called the {\em center} of the fan, and $2n$ other vertices $v_1,\dots, v_{2n}$ such that for $i = 1,\dots,n$, $vv_{2i-1}v_{2i}$ is a triangle. Each of the $n$ edges $v_{2i-1}v_{2i}$ is called a {\em blade} of the fan.

For the Ramsey numbers of $R(C_m,C_n)$, it has been studied and completely determined in Bondy and Erd\H{o}s \cite{b-e}, Faudree and Schelp \cite{f-s}, and Rosta \cite{ros}. 		
Ramsey numbers of fans $R(F_m, F_n)$ have been studied, both in the diagonal case (when $m = n$) and the off-diagonal case.
 For results in the off-diagonal case, see \cite{Li-ro-1996, Lin-li-2009, Lin-li-2010, zhang-bro-2015}. In particular, Lin, Li and Dong \cite{Lin-li-2010} showed that $R(F_m,F_n)=4n+1$ for each fixed $m\ge1$ and large $n$.
 Recently,  Chen, Yu and Zhao \cite{chen-yu} improve the bounds for $R(F_n,F_n)$ significantly and obtain  that
 $$\frac{9n}{2}-5\leq R(F_n,F_n)\leq \frac{11n}{2}+6,$$
 and Dvo\v{r}\'{a}k and Metrebian \cite{Metr-2021} make a further improvement on the upper bound by decreasing the coefficient of the main term from $5.5$ to about $5.167$.

%For the Ramsey numbers of $R(C_m,C_n)$, it has been studied and completely determined in Bondy and Erd\H{o}s \cite{b-e}, Faudree and Schelp \cite{f-s}, and Rosta \cite{ros}. 	

The Ramsey numbers of cycles versus fans also attracted much of attention. For instance, Li and Rousseau \cite{Li-ro-1996} obtained that $R(C_3,F_n)=4n+1$ for all $n\ge2$, one can see also in Bollobas \cite[Theorem 13 in Ch. 6]{bol}. Generally, for fixed $m$ and large $n$, Liu and Li \cite{Liu-li-2016} showed that $R(C_{2m+1},F_n)=4n+1$. Shi \cite{Shi-2010} considered the case when the order of cycle is much larger than that of fan,  in particular, the author showed that $R(C_n,F_m)=2n-1$ holds for all $n>3m$.
For more Ramsey numbers involving fans, we refer the reader to \cite{ bre-2017, Rad-2014, Salman-2006,  zhang-2015}, etc.
%\cite{b-e, f-s, ros, nik-2004, nik-2009, }	
	
In this paper, we are concerned with the asymptotic behavior of the Ramsey number $R(C_{2\lfloor an\rfloor}, F_n)$ when  $n$  is large and $a\geq 1/2 $ is fixed.
\begin{theorem}\label{main theorem}
For all sufficiently large $n$,
\[
R(C_{2\lfloor an\rfloor}, F_n)= \left\{ \begin{array}{ll}
(2+2a+o(1))n & \textrm{if $1/2\leq a< 1$,}\\
(4a+o(1))n & \textrm{if $ a\geq 1 $.}
\end{array} \right.
\]
\end{theorem}

The following corollary is immediate.
\begin{corollary} We have
$R(C_{2n}, F_n)=(4+o(1))n$, and $R(C_{n}, F_n)=(3+o(1))n$ for all sufficiently large even integer $n$.
\end{corollary}

\section{Preliminaries}\label{chap2}

Throughout this paper, all graphs are finite and simple. Let $G=G(V,E)$ be such a graph.
%If $e = \{u, v\}\in E$ (in short, $e = uv$), then $u$ is called adjacent to $v$, and $u$ and $v$ are called neighbors.
For a vertex $v\in V$, let $N_G(v)$ denote the neighborhood of $v$ in $G$,
and $\deg_G(v)=|N_G(v)|$ is the degree of a vertex $v\in V$.
We denote by $\delta(G)$  and $\triangle(G)$ the minimum and maximum degrees of the vertices of $G$.
For a vertex $v\in V$ and $U\subset V$, we write $N_G(v,U)$ for the neighbors of $v$ in $U$ in graph $G$ and denote $\deg(v,U)=|N_G(v,U)|$.
For a vertex set $X\subset V$ and $U\subset V\setminus X$, we write $N_G(X,U)$ for all neighbors of $X$ in $U$ in graph $G$.
In particular, we write $N_G(X)$ for all neighbors of $X$ in  $V\setminus X$.

For disjoint vertex sets $A, B\subseteq V$, let $e_G(A,B)$ denote the number of edges of $G$ with one endpoint in $A$ and the other in $B$, and the density between $A$ and $B$ is
	\[d_G(A,B)=\dfrac{e_G(A,B)}{|A| |B|}.\]
	We always delete the subscript when there is no confusion.
	\begin{definition}[$\epsilon$-regular]\label{regular}
		For $\epsilon>0$ and $d\leq 1$, a pair $(A,B)$ is $\epsilon$-regular if for all $X\subseteq A$ and $Y\subseteq B$ with $|X|>\epsilon|A|$ and $|Y|>\epsilon|B|$ we have $|d(X,Y)-d(A,B)|<\epsilon$.
	\end{definition}
	\begin{definition}[($\epsilon,d$)-regular]
		A pair $(A,B)$ is said to be $(\epsilon,d)$-regular if it is $\epsilon$-regular and $d(A,B)\geq d$.
	\end{definition}
%	\begin{definition}[$(\epsilon,d)$-super-regular]
%		A pair $(A,B)$ is said to be $(\epsilon,d)$-super-regular if it is $\epsilon$-regular and $\deg(u,B)>d|B|$ for all $u\in A$ and $\deg(v,A)>d|A|$ for all $v\in B$.
%	\end{definition}
	
The following property is well-known, see e.g. \cite{kom}.	
%Let us have a property that any regular pair has a large subgraph which is super-regular, and we include a proof for completeness.
\begin{fact}\label{most degree large}
Let $(A,B)$ be an $\epsilon$-regular pair with density $d$. Then for any  $Y\subset B$, $|Y|>\epsilon |B|$ we have $$\#\left\{x\in A: \deg(x,Y)\leq(d-\epsilon)|Y| \right\}\leq \epsilon |A|.$$
\end{fact}

% \begin{fact}\label{hdt}
% Let $(A,B)$ be an $\epsilon$-regular pair in $G$ with density $d$, and let $X\subseteq A$ and $Y\subseteq B$ with $|X|\geq \alpha|A|$ and $|Y|\geq \alpha|B|$ for some $\alpha>\epsilon$. Then $(X,Y)$ is $\epsilon'$-regular pair with  $\epsilon'=\max\{\epsilon/\alpha, 2\epsilon\}$, and   for its density $d'$ we have $|d'-d|<\epsilon$.
%\end{fact}

%\begin{fact}\label{localregular}
%For $0<\epsilon<1/2$ and $d\leq1$, if $(A,B)$ is $(\epsilon,d)$-regular with $|A|=|B|=m$, then there exist $A_1\subseteq A$ and $B_1\subseteq B$ with $|A_1|=|B_1|\ge(1-\epsilon)m$ such that $(A_1,B_1)$ is $(2\epsilon,d-2\epsilon)$-super-regular.
%\end{fact}

\medskip
		
In this paper, we will use the following regularity lemma.

\begin{lemma}[Szemer\'{e}di \cite{regular-lemma}]\label{regular lemma}
For every $\epsilon> 0$ and integer $t_0\geq 1$, there exists $T_0=T_0(\epsilon,t_0)$
such that,   for every graph $G$  of large order $n$,
there exists a partition $V (G)= \cup_{i=0}^t{V_i}$ satisfying $t_0\le t\le T_0$ and

\smallskip
		
$(i)$ $\left| {{V_0}} \right| < \epsilon n$, $\left| {{V_1}} \right| = \left| {{V_2}} \right| =  \ldots  = \left| {{V_t}} \right|$;
		
\smallskip
		
$(ii)$  all but at most $\epsilon{t^2}$ pairs $(V_i,V_j)$, $1 \leq i\neq j\leq t$, are $\epsilon$-regular.
\end{lemma}

%The next lemma by Benevides and Skokan \cite{Benevides-Skokan-2009} is a slightly stronger version compared to the original one established by {\L}uczak \cite[Claim 3]{lucazk-1999}.
%	
%\begin{lemma}[Benevides and Skokan \cite{Benevides-Skokan-2009}]\label{long-path-lemma}
%For every $0<\beta_0<1$, there exists an $n_0$ such that for every $n>n_0$ the following holds: If $(V_1,V_2)$ is $\epsilon$-regular with $|V_1|=|V_2|=n$ and density at least $\beta_0/4$  for some $\epsilon$ satisfying $0<\epsilon<\beta_0/100$, then for every $\ell, 1\leq \ell \leq n-5\epsilon n/\beta_0$, and for every pair of vertices $v'\in V_1$, $v''\in V_2$ satisfying $\deg(v',V_2)$, $\deg(v'', V_1)\geq \beta_0 n/5$, $G$ contains a path of length $2\ell+1$ connecting $v'$ and $v''$.
%\end{lemma}

For a graph $G$, denote by $\nu(G)$ the size of the largest matching of $G$. Let us recall the following classical result in graph theory due to Hall,   see, e.g., \cite{b-m,west}.
\begin{lemma}[Hall \cite{b-m}]\label{hall}
Let $G$ be a bipartite graph on parts $X$ and $Y$. For any non-negative
integer $d$, $\nu(G) \geq |X|-d$  if and only if $|N(S)| \geq |S|-d $ for every
$S \subseteq X$.
\end{lemma}

 For a matching $M \subseteq E$, we call all vertices
which are not incident to any edge in $M$ the unmatched vertices in $M$. Furthermore, we
denote by $q(G\setminus S)$ the number of odd components in $G \setminus S$. We use a generalization of
Tutte's Theorem  in our proof.
\begin{lemma}[Berge \cite{Berge-1958}]\label{ber-com}
Let $G = (V, E)$ be a graph. For any set $S\subseteq V$ and any
matching $M$, the number of unmatched vertices in $M$ is at least
$q (G\setminus S)-|S|$. Moreover, there exists a set $S\subseteq V$ such that every maximum matching of $G$ misses exactly $q (G\setminus S)-|S|$
vertices.
\end{lemma}

For a graph $G$, we use $g(G)$ and $c(G)$ to denote its \emph{girth} and \emph{circumference}, i.e., the length of a shortest cycle and a longest cycle of $G$.
We say a graph is {\it 2-connected} if it remains connected after the deletion of any  vertex.

%The following classical result is due to Dirac \cite{di}.
\begin{lemma}[Dirac \cite{di}]\label{1-thm-3}
Let $G$ be a 2-connected graph of order $n$ with minimum degree $\delta=\delta(G)$. Then $c(G)\geq \min\{2\delta,n\}$.
\end{lemma}

A graph is called \emph{weakly pancyclic} if it contains cycles of every length between its  girth and its circumference.
A graph is \emph{pancyclic} if it is weakly pancyclic with girth $3$ and circumference $n=|V(G)|$.
In particular, if $\delta=\delta(G)\ge n/2$, then $c(G)=n$. This is a well-known result for a graph being hamiltonian.
For the special case of $\delta \geq n/2$,
the following result  tells us more about the structure of a graph.

\begin{lemma}[Bondy \cite{Bondy}]\label{1-thm-1}
If a graph $G$ with $n$ vertices satisfies  $\delta(G)\geq n/2$, then $G$ is pancyclic unless $n=2r$ and $G=K_{r,r}$.
\end{lemma}

Let $nK_2$ denote a matching of size $n$, i.e., $n$ pairwise disjoint edges, and let $S_t$ be a star with $t$ edges.
%The following Result obtained by Gy\'{a}rf\'{a}s and S\'{a}rk\"{o}zy      \cite{gya-sar-2012} in 2012.

\begin{lemma}[Gy\'{a}rf\'{a}s and S\'{a}rk\"{o}zy \cite{gya-sar-2012}]\label{star-match}
Suppose that $n_1 \geq n_2 \geq 1$ and $k \geq 1$. Then
\begin{equation*}
R(S_k, n_1K_2, n_2K_2)=\left\{
\begin{aligned}
&2n_1+n_2-1~&\text{if}~ k\leq n_1,\\
& n_1+n_2-1+t~&\text{if}~ k\geq n_1.
\end{aligned}
\right.
\end{equation*}
\end{lemma}

%The following result due to Balogh, Kostochka, Lavrov and   Liu \cite{balogh-2019} confirms a conjecture of Schelp \cite{conj-schelp}.
%\begin{lemma}[Balogh et al. \cite{balogh-2019}]\label{balogh-conj}
%Let $G$ be a graph on $3n -1$ vertices with minimum degree at least $(3|V(G)|-1)/4$. If $n$ is sufficiently large, then $G\rightarrow(P_{2n},P_{2n})$.
%\end{lemma}

%Erd\H{o}s and Rado remarked that any 2-colored complete graph contains a monochromatic spanning tree, see \cite{gy}.
We will apply the following result to get a large monochromatic component for every $2$-coloring of the edges of graph $G$ with large minimum degree.

\begin{lemma}[Gy\'{a}rf\'{a}s and S\'{a}rk\"{o}zy \cite{gya-sar-2012}]\label{gya-lemma}
For any $2$-color of edges of a graph $G$ with minimum degree $\delta(G)\geq\frac{3|V(G)|}{4}$,
there is a monochromatic component of order larger than $\delta(G)$. This estimate is sharp.
\end{lemma}

We also need the following result, which  states that a bipartite graph with high density always contains a large matchings.

\begin{lemma}[Figaj  and {\L}uczak \cite{fl}]\label{bi-match}
Let $G=(V,E)$ be a bipartite graph with bipartition $\{V_1,V_2\}$,
$|V_1|\geq|V_2|$, and at least $(1-\epsilon)|V_1||V_2|$ edges, for some $0<\epsilon<0.01$.  Then, there is a component in $G$ of at least $(1-3\epsilon)(|V_1|+|V_2|)$ vertices which contains a matching of cardinality at least $(1-3\epsilon)|V_2|$.
\end{lemma}

\section{Proof of Theorem \ref{main theorem} }

For a graph $G=G(V,E)$, if $S \subseteq V$, then $G - S$  denotes the subgraph
of $G$ induced by $V (G)\setminus S$.
 For any subset $A \subseteq V$, we use $G[A]$ to denote the subgraph induced by the vertex set $A$ in $G$.
For two subsets $A \subseteq V$ and $B \subseteq V$, we use $G[A,B]$ to denote the subgraph induced by all edges between $A$ and $B$ in $G$.

For a 2-edge colored graph $G$, we use $G^r$ (or $G^b$) to denote the subgraph of $G$ formed by all red (or blue) edges of $G$.
 If $A\subseteq V$, then $A^r$ is defined to be $G^r[A]$ and $A^b$  is defined to be $G^b[A]$.
 For any subset $A\subseteq V$, we also use $\delta^r(A) $ (or $\delta^b(A) $) to denote $\delta(A^r)$ (or $\delta(A^b)$) for convenience.

In the following, we always omit the floors and ceilings when there is no affection on our argument.

\bigskip
{\bf Part (I)  \ $1/2\le a<1$}

\bigskip
Let $\overline{G}$  denote the complement graph of $G$.
Note that  the graph $K_{2 \lfloor an\rfloor-1}\cup {K_{n-2}}\cup {K_{n-2}}$
 contains  no cycle $C_{2 \lfloor an\rfloor}$ and its complement contains no $F_n$, so we have $R(C_{2 \lfloor an\rfloor}, F_n)\geq 2 \lfloor an\rfloor +2n-4$ for $1/2\leq a< 1$.

 It remains to show the upper bound. Let $N=(2a+2+\gamma)n$, where $1/2\leq a< 1$ and $0<\gamma<1/10$ is sufficiently small, we will show $R(C_{2an},F_n)\leq N$ for all large $n$, i.e.,
any red/blue edge coloring of $K_N$ yields either a red $C_{2an}$ or a blue  $F_n$.
Suppose to the contrary that for large   $n$, there exists a coloring that contains neither a red $C_{2an}$ nor a blue $F_n$. We aim to find a contradiction.

Consider a 2-edge coloring of $G=K_N$ defined on $V$. Set
\begin{align}\label{eta-ep}
\beta=\min\left\{\frac{\gamma}{100}, \frac{1-a}{30(a+1)}\right\},
\;\;\text{and}\;\;\epsilon =\frac{\beta^2}{10^4}.
\end{align}
We apply the regularity lemma (Lemma \ref{regular lemma}) with $\epsilon$ and sufficiently large $t_0$ to obtain
\begin{align}\label{c-2}
T_0=T_0(\epsilon,t_0)=\min\left\{2t_0, \frac{5}{4\epsilon}\right\}
\end{align}
such that there exists a partition  $V=V_{0}\cup V_{1}\cup\dots \cup V_t$ satisfying $t_0\le t\le T_0$ and
	(i) $\left| {{V_0}} \right| < \epsilon N$, $\left| {{V_1}} \right| = \left| {{V_2}} \right| =  \ldots  = \left| {{V_t}} \right|$;
	(ii) all but at most $\epsilon t^ 2$ pairs $(V_i,V_j)$, $1 \le i \neq j \le t$, are $\epsilon$-regular for $G^r$ and $G^b$.
We construct the reduced graph $H$ with vertex set $\{v_1, v_2, \dots, v_t\}$ and the edge set formed by pairs $\{v_i,v_j\}$ for which $(V_i,V_j)$ is $\epsilon$-regular with respect to $G^r$ and $G^b$.
%In which $v_i$ and $v_j$ are non-adjacent in $H$ if the pairs $(V_i, V_j)$ is not  $\epsilon$-regular.
Thus we obtain a bijection $f:{v_i} \to {V_i}$ between the vertices of $H$ and the clusters of the partition.

Color an edge $v_{i}v_{j}$ red if the density of the red edges between $V_i$ and $V_j$ is at least $\beta$, and blue otherwise. Let $H^r$ and $H^b$ be the subgraphs induced by all red edges and blue edges of $H$, respectively.
Since there are at most $\epsilon t^2$ edges that are uncolored in $H$, by deleting at most $\sqrt{\epsilon}t$ vertices, we may assume that each vertex is adjacent to at most $\sqrt{\epsilon}t$ non neighbors. In what follows, when referring to the reduced graph $H$, we will assume that these vertices have been removed.
	
%In the following, we shall omit all floors and ceilings since it will not affect our argument.

A connected matching in a graph $G$ is a matching $M$ such that all edges of $M$ are in the same connected component of $G$.

\begin{claim}\label{long-match-1}
$H^r$ contains no  connected matching of size more than $(\frac{a}{2a+2}-0.15\beta)t$.
\end{claim}
\noindent{\bf Proof.}
On the contrary, suppose that $H^r$ contains a connected matching $M$
on at least $2k=(a/(a+1)-0.3\beta)t$ vertices. Let $F$ be a minimal connected red subgraph containing $M$. We may assume that $M=\{v_{1}v_2,v_3v_4,\dots,v_{2k-1}v_{2k}\}$ and $f(v_{2i-1})=V_{2i-1}$ for $1\leq i\leq k$. Clearly, $F$ is a tree.
Consider a closed walk $W=V_1V_2\cdots V_{2k-1}V_{2k}\cdots V_1$ that contains all edges of $M$.
%Since $F$ is a tree, $W$  must be of even length.
By applying a similar argument as in Figaj and {\L}uczak \cite{fl} we can obtain a red cycle of length $2an$, contradicting our assumption that $G^r$ contains no such cycle.
%This completes the proof of Claim \ref{long-match-1}.
\hfill$\Box$

%Since $G^b$ contains no $F_n$, this implies that the following claim.	

\begin{claim}\label{large-fan-1}
$H^b$ contains no fan with at least $k=(\frac{1}{2a+2}-0.05\beta)t$ blades.
\end{claim}

\noindent{\bf Proof}.
If not,  $H^b$ contains  a fan on $2k+1$ vertices.
Suppose that $v_a$ is the center of such fan with $k$ blades, say $v_1v_2,\dots,v_{2k-1}v_{2k}$.
By relabelling the vertices if needed, we may assume that $f(v_a)=V_a$ and $f(v_i)=V_i$ for $i\in [2k]=\{1,2,\dots,2k\}$.
%We will prove that there exists a  copy of $F_n$ in $G^b$, contradicting the assumption that $G^b$ contains no $F_n$.

Note that an edge  $v_iv_j$ in $H$ is blue if and only if the density $d_{G^b}(V_i,V_j)\geq 1-\beta$,
by Fact \ref{most degree large}, all but at most $2k\epsilon$ vertices of $V_a$ has degree at least $(1-\beta-\epsilon)|V_i|$  in each $V_i$ for $1\leq i\leq 2k$.
Since $2k\epsilon<1$ from (\ref{c-2}), we can choose a vertex $u\in V_a$ such that $u$ has at least $(1-\beta-\epsilon)|V_i|$ neighbors in each $V_i$ for $1\leq i\leq 2k$.
Let $V_i'=N_{G^b}(u)\cap V_i$  for $1\leq i\leq 2k$.
Therefore,
\begin{align}\label{ineq-6}
|V_i'|=|N_{G^b}(u)\cap V_i|\geq (1-\beta-\epsilon)|V_i|>\epsilon|V_i|
\end{align}
for every $1\leq i\leq 2k$.
Moreover, we have $e_{G^b}(V_{2i-1}',V_{2i}')\geq(1-\beta-\epsilon)|V_{2i-1}'||V_{2i}'|$
for $1\leq i \leq k$ since $(V_{2i-1},V_{2i})$ are $(\epsilon,1-\beta)$-regular in $G^b$.
Hence, by Lemma \ref{bi-match}, the graph $G^b[V_{2i-1},V_{2i}]$ contains a matching of cardinality at least
$
\left( 1-3(\beta+\epsilon) \right)|V_{2i-1}'|.
$
Let $S =  \cup_{i = 1}^{2k} {{V_i'}}$.
Therefore, $G^b[S] $ contains  a  matching  of cardinality at least
\begin{align*}
k ({1 - 3\left( {\beta  + \epsilon } \right)})\left| {{V}_i'} \right|
&\overset{(\ref{ineq-6})}{\ge} k ( {1 - 3\left( {\beta  + \epsilon } \right)} )(1 - \left( {\beta  + \epsilon } \right))\left| {{V_i}} \right|\\
&\ge k( {1 - 4\left( {\beta  + \epsilon } \right)} )\frac{{(1 - \epsilon )(2a+2 + \gamma )n}}{t}\\
& \overset{(\ref{eta-ep})}{>} \left ( {\frac{1}{2a+2} - 0.05\beta } \right)\left( {1 - 5\beta } \right)\left( {2a+2+ \gamma } \right)n \\
& > \left( {\frac{1}{2a+2}- 2\beta } \right)\left( {2a+2 + \gamma } \right)n \overset{(\ref{eta-ep})}{>} \left( {1 + 2\beta} \right)n,
\end{align*}
yielding a  blue $F_n$ with center $u$ in $G^b$, a contradiction.
\hfill$\Box$

\begin{claim}\label{cla-match-1}
$\deg_{H^r}(v)\geq \frac{a}{2a+2}t+1$ for every $v\in V(H)$.
\end{claim}

\noindent{\bf Proof.}
On the contrary, without loss of generality, suppose that $H^r$ contains a vertex $v_c $ such that $\deg_{H^r}(v_c)\leq \frac{a}{2a+2}t$.
Since $\deg_H(v)\geq(1-2\sqrt{\epsilon})t-1$ for each $v\in V(H)$,
we have that $\deg_{H^b}(v_c)\geq (\frac{a+2}{2a+2}-2\sqrt{\epsilon})t-1$.
Denote $H_1=H[N_{H^b}(v_c)]$, i.e., the subgraph of $H$ induced by the neighborhood of $v_c$ in $H^b$.
Without loss of generality,  we may assume that
\begin{align}\label{ineq-15}
|V(H_1)|=\left(\frac{a+2}{2a+2}-2\sqrt{\epsilon}\right)t-1.
\end{align}
Note that every vertex in $H$ has at most $\sqrt{\epsilon}t$ non-neighbors.
Let $C_1$ be the vertex set of  a largest monochromatic component in $H_1$. Here and in what follows, we also use $C_1$ ($C_i$) to denote its vertex set of the component $C_1$ ($C_i$).
%Observe that the largest  matching in $C_1$ is automatically connected.
From Lemma  \ref{gya-lemma},
\begin{align}\label{ineq-1}
|C_1|>\delta(H_1)\ge|V(H_1)|-1-\sqrt{\epsilon}t\overset{(\ref{ineq-15})}{= }
\left(\frac{a+2}{2a+2}-3\sqrt{\epsilon}\right)t-2.
\end{align}

Suppose first that $C_1$ is a red component.
We  apply Lemma \ref{star-match} to the subgraph $C_1$  with $k=\sqrt{\epsilon}t+1$, $n_1=(\frac{1}{2a+2}-0.05\beta)t$, and $n_2=(\frac{a}{2a+2}-0.15\beta)t$ to obtain that
\[
R(S_k, n_1K_2, n_2K_2)=2n_1+n_2-1=\left(\frac{a+2}{2a+2}-0.25\beta\right)t-1\overset{(\ref{ineq-1}), (\ref{eta-ep})}{<}|C_1|.
\]
Thus, by noting that every vertex in $H$ has at most $\sqrt{\epsilon}t$ non-neighbors, we can get either
 a red matching of size $(\frac{a}{2a+2}-0.15\beta)t$ or a blue matching of size $(\frac{1}{2a+2}-0.05\beta)t$.
The first case contradicts Claim \ref{long-match-1} since the red matching in $C_1$ is clearly connected.
For the second, we obtain a blue fan with center $v_c$ in $H_1^b$ with at least
$(\frac{1}{2a+2}-0.05\beta)t$ blades, which contradicts Claim \ref{large-fan-1}.

In the following, we assume that $C_1$ is a largest blue component.
Denote $U=V(H_1)\setminus V(C_1).$
Recall that every vertex in $H_1$ has at most $\sqrt{\epsilon}t$ non-neighbors.
If $U\neq{\O}$, then  $U$ is completely covered
by a red component $C_r$ due to the minimum degree condition of $H_1$.
Note that there are no  edges of $H_1$ in the bipartite graph $H_1[C_1\setminus C_r, C_r\setminus C_1]$,
we have that $|C_1\setminus C_r|\leq\sqrt{\epsilon}t$ since $C_r\setminus C_1=U\neq {\O}$. %By (\ref{ineq-1}), so we have $|C_2\setminus C_1|<\sqrt{\epsilon}t$.
Thus we have
\[
|C_1\cap C_r|=|C_1|-|C_1\setminus C_r|\overset{(\ref{ineq-1})}{\geq}\left(\frac{a+2}{2a+2}-
4\sqrt{\epsilon}\right)t-1\overset{(\ref{eta-ep})}{>}\left(\frac{a+2}{2a+2}-
0.35\beta\right)t-1.
\]
We  apply Lemma \ref{star-match} to the subgraph induced by $C_1\cap C_r$ in $H_1$ to conclude that $H_1[C_1\cap C_r]$
contains a red connected matching of size at least $(\frac{a}{2a+2}-0.15\beta)t$
or a blue matching of size at least $(\frac{1}{2a+2}-0.1\beta)t$ which together with $v_c$ yield a blue fan with more than $(\frac{1}{2a+2}-0.1\beta)t$ blades, contradicting Claim \ref{long-match-1} or Claim \ref{large-fan-1}.

Now we assume $U={\O}$. So we have $C_1=V(H_1)$ and hence
\begin{align}\label{ineq-12}
|C_1|=|V(H_1)|\overset{(\ref{ineq-15})}{= } \left(\frac{a+2}{2a+2}-2\sqrt{\epsilon}\right)t-1.
\end{align}
Without loss of generality, we  define $C_2$ as a largest red component in $H_1$.
From Claim \ref{large-fan-1}, we know that the largest blue matching has
size $m<(\frac{1}{2a+2}-0.05\beta)t$ in $H_1$.
Applying Lemma \ref{ber-com}  to the subgraph induced by all blue edges in $H_1$, we can find a subset $S\subset V(H_1)$ such that the number of odd components
\begin{align}\label{ineq-7}
q(V(H_1)\setminus S)= |S|+|V(H_1)|-2m
\nonumber  &\geq\left(\frac{a+2}{2a+2}-2\sqrt{\epsilon}\right)t-2\left(\frac{1}{2a+2}-0.05\beta\right)t\\  &\overset{(\ref{eta-ep})}{>}\left(\frac{a}{2a+2}+6\sqrt{\epsilon}\right)t.
\end{align}
Clearly, $q(V(H_1)\setminus S)+|S|\leq |V(H_1)|$, which implies that
$|S|\leq m\leq (\frac{1}{2a+2}-0.05\beta)t-1$.
Let $R$ be the red subgraph of $H_1$ whose
vertex set is $V(H_1)\setminus S$ and edge set consists of all red edges between
blue components of $V(H_1)\setminus S$ in $H_1$.
It is clear that
\begin{align}\label{ineq-8}
|V(R)|\geq |V(H_1)|-|S|\geq \left(\frac{a+2}{2a+2}
-2\sqrt{\epsilon}\right)t-\left(\frac{1}{2a+2}-0.05\beta\right)t\overset{(\ref{eta-ep})}{>}\left(\frac{1}{2}+2\sqrt{\epsilon}\right)t.
\end{align}

We will show that $R$ is connected. Otherwise, $V(R)$ can be partitioned into two non-empty
sets $A$ and $B$ such that there are no red edges between $A$ and $B$. Without loss of generality, suppose $|A|\ge|B|$. Then $|A|> (1/4+\sqrt{\epsilon})t$.
If $A$ intersects each of these $q(V(H_1)\setminus S)$  blue odd components in $V(H_1)\setminus S$,
then any vertex $v \in B$ is non-adjacent to all vertices in the intersecting set of $A$ and
 those blue components not containing $v$. Thus, any vertex $v \in B$ is non-adjacent to at least
\[
q(V(H_1)\setminus S)-1\overset{(\ref{ineq-7})}{>}\left(\frac{a}{2a+2}+6\sqrt{\epsilon}\right)t-1
\]
vertices in $H_1$.
On the other hand, if $A$ does not intersect some blue odd component, then any vertex $v$ in this component is non-adjacent to any vertex of $A$.
Therefore, in both cases we can find a vertex $v$ that is non-adjacent to at least $(\frac{a}{2a+2}+6\sqrt{\epsilon})t-1$ vertices, which clearly contradicts the fact that $\delta(H_1)\geq |V(H_1)|-1-\sqrt{\epsilon}t$. Thus, $R$ is connected as desired.
Since $C_2$ is the largest  red component, it follows from (\ref{ineq-8}) that
\begin{align}\label{size-C2}
|C_2|\geq |V(R)|>(1/2+2\sqrt{\epsilon})t.
\end{align}

Let $p=|C_1\setminus C_2|$.
Then $p\geq 23\sqrt{\epsilon}t$. Otherwise, $|C_2|>(\frac{a+2}{2a+2}-0.25\beta)t-1$ from (\ref{ineq-12}), by a similar argument as above by applying Lemma \ref{star-match}, we can get either a red connected matching of size at least $(\frac{a}{2a+2}-0.15\beta)t$, or a blue matching of size at least $(\frac{1}{2a+2}-0.05\beta)t$
which together with $v_c$ yield a blue fan with more than $(\frac{1}{2a+2}-0.05\beta)t$ blades. This again leads to a contradiction from Claim \ref{long-match-1} or Claim \ref{large-fan-1}.

We first suppose that $23\sqrt{\epsilon}t \leq p<(\frac{1-a}{a+1}-0.05\beta)t$.
Thus
\begin{align}\label{ineq-2}
|C_2|&= |C_1|-p\overset{(\ref{ineq-12})}{=}
\left(\frac{a+2}{2a+2}-2\sqrt{\epsilon}\right)t-1-p.
\end{align}
We apply Lemma \ref{star-match} to the subgraph induced by vertex set $C_2$ in $H_1$ with $k=\sqrt{\epsilon}t+1$, $n_1=(\frac{1}{2a+2}-0.05\beta)t-p/2$ and $n_2=(\frac{a}{2a+2}-0.15\beta)t$. Note that $n_1> n_2$ and
\[
R(S_k, n_1K_2, n_2K_2)=\left(\frac{a+2}{2a+2}-0.25\beta\right)t-p-1\overset{(\ref{ineq-2})}{<}|C_2|,
\]
so there exists a blue matching $M_1$ of size $n_1=(\frac{1}{2a+2}-0.05\beta)t-p/2$ since otherwise a red connected matching of size at least $n_2$ will lead to a contradiction from Claim \ref{long-match-1}.

Note that $\delta(H_1)=\delta(H_1[C_1])\geq |C_1|-1-\sqrt{\epsilon}t$ and
all (but at most $\epsilon t^2$)  edges between $C_1\setminus C_2$ and $C_2\setminus V(M_1)$ are blue.
For any subset $S\subseteq C_2\setminus V(M_1)$,
if $$|C_2\setminus V(M_1)|<|C_1\setminus C_2|=p\leq\left(\frac{1-a}{a+1}-0.05\beta\right)t,$$
then the total number of blue neighbors of $S$ in $C_1\setminus C_2$ satisfies that
\[
|N_{H^b_1}(S, C_1\setminus C_2)|\geq  |C_1\setminus C_2|-\sqrt{\epsilon}t>|C_2\setminus V(M_1)|-\sqrt{\epsilon}t.
\]
Recall that $|C_2|=|C_1|-p$ and $n_1=(\frac{1}{2a+2}-0.05\beta)t-p/2$.
Hence, by Lemma \ref{hall}, the bipartite graph $H^b_1[C_1\setminus C_2, C_2\setminus V(M_1)]$ contains
a blue matching  of size at least
\begin{align*}
|C_2\setminus V(M_1)|-\sqrt{\epsilon}t&=|C_2|-2n_1-\sqrt{\epsilon}t=|C_1|-\left(\frac{1}{a+1}-0.1\beta\right)t-\sqrt{\epsilon}t\\
&\overset{(\ref{ineq-12})}{=}
\left(\frac{a}{2a+2}+0.1\beta-3\sqrt{\epsilon}\right)t-1
\\&\overset{(\ref{eta-ep})}{>}\left(\frac{a}{2a+2}+0.06\beta\right)t.
\end{align*}
This matching together with $M_1$ yield a blue  matching of size at least
\begin{align*}
n_1+\left(\frac{a}{2a+2}+0.06\beta\right)t=(1/2+0.01\beta)t-p/2>\left(\frac{1}{2a+2}-0.05\beta\right)t
\end{align*}
since $p<(\frac{1-a}{a+1}-0.05\beta)t$.
If $|C_2\setminus V(M_1)|\geq p=|C_2|$, then Lemma \ref{hall} again implies the bipartite graph
$H^b_1[C_1\setminus C_2, C_2\setminus V(M_1)]$ contains
a blue matching  of size at least
\[
|C_1\setminus C_2|-\sqrt{\epsilon}t=p-\sqrt{\epsilon}t,
\]
which together with $M_1$ yield a blue  matching of size at least $n_1+p-\sqrt{\epsilon}t>(\frac{1}{2a+2}-0.05\beta)t$
in $H_1^{b}$ (also in $H^{b}$).
Therefore, for either case, we can get a blue fan with center $v_c$ and at least $(\frac{1}{2a+2}-0.05\beta)t$ blades, which contradicts Claim \ref{large-fan-1}.

Now we assume that $p\geq(\frac{1-a}{a+1}-0.05\beta)t$.
Recall that $C_1=V(H_1)$ and $|C_2|>(1/2+2\sqrt{\epsilon})t$ from (\ref{size-C2}),
so we can upper bound $p$ as that
 \begin{align*}
 p=|C_1\setminus C_2|=|C_1|-|C_2| \overset{(\ref{ineq-12})}{<}\left(\frac{1}{2a+2}-4\sqrt{\epsilon}\right)t
 \overset{(\ref{eta-ep})}{=}\left(\frac{1}{2a+2}-0.04\beta\right)t.
 \end{align*}
Thus we have $(\frac{1-a}{a+1}-0.05\beta)t\leq p<(\frac{1}{2a+2}-0.04\beta)t$.
We apply Lemma \ref{star-match} to the subgraph induced by $C_2$ in $H_1$ with $k=\sqrt{\epsilon}t+1$, $n_1=(\frac{a}{2a+2}+0.01\beta)t$ and $n_2=(\frac{1}{2a+2}-0.04\beta)t-p$.
 Note that $n_1\geq n_2$ and $1/2\leq a<1$, hence we have that
\[
R(S_k, n_1K_2, n_2K_2)=\left(\frac{2a+1}{2a+2}-0.02\beta\right)t-p-1\overset{ (\ref{eta-ep}),(\ref{ineq-12})}{\le} |C_1|-p=
|C_2|.
\]
Thus there is a blue matching $M_2$ of size $n_2=(\frac{1}{2a+2}-0.04\beta)t-p$ since otherwise a red
connected matching of size at least  $n_1=(\frac{a}{2a+2}+0.01\beta)t$
will again lead to a contradiction from Claim \ref{long-match-1}.

It is clear that $|C_2\setminus V(M_2)|>p+\sqrt{\epsilon}t$.
Recall  that $\delta(H_1)=\delta(H_1[C_1])\geq |C_1|-1-\sqrt{\epsilon}t$ and
all edges between $C_1\setminus C_2$ and $C_2\setminus V(M_2)$ are blue.
Thus, by Lemma \ref{hall}, the bipartite graph $H^b_1[C_1\setminus C_2, C_2\setminus V(M_2)]$ contains
a blue matching  of size $p$, which together with $M_2$ yield a blue  matching of size  $n_2+p=(\frac{1}{2a+2}-0.04\beta)t$ in $H_1^{b}$.
Therefore, we can get a blue fan with center $v_c$ and at least $(\frac{1}{2a+2}-0.05\beta)t$ blades, which contradicts Claim \ref{large-fan-1}.

This completes the proof of Claim \ref{cla-match-1}.
\hfill$\Box$

\begin{claim}\label{2-connect}
$H^r$ is  2-connected.
\end{claim}

\noindent{\bf Proof. }
Suppose $H^r$ is not 2-connected. Then there exists a $K\subseteq V(H)$ with $|K|\leq 1$ such that $(H-K)^r$ is disconnected. Let
$C_1,C_2,\dots,C_\tau$ be the  vertex sets of red components of $(H-K)^r$, such that
$|C_1|\geq |C_2|\dots \geq|C_\tau|$. Note that $\delta^r(H-K)\geq \delta (H^r)-|K|\geq \frac{a}{2a+2}t$ by Claim \ref{cla-match-1},
thus for $1\le i\le \tau$,
\begin{align}\label{ineq-9}
|C_i|\ge\frac{a}{2a+2}t+1, \;\; \text{and}\;\;\delta^r(C_i)\geq \delta(H^r)-|K|\geq \frac{a}{2a+2}t
\end{align}
for each $i\in[\tau]$.
Therefore,  there are at most five components in $(H-K)^r$, i.e., $\tau\leq 5$, by noting $1/2\leq a<1$.
%Note that  all edges  joining distinct components are all blue.
The proof  is divided into four cases according to the number of red components.

\medskip

\noindent\textbf{Case A:} \ $\tau=5$

\medskip

%Clearly, all edges joining distinct $C_i$, $C_j$ are  blue and $| C_i|\geq \frac{a}{2a+2}t+1$ for every $1\leq i\leq 5$.
Since each vertex $v\in V(H)$ has at most $\sqrt{\epsilon}t$ non-neighbors,
it follows from Lemma \ref{hall} that the bipartite graph $H^b[C_1, C_2]$ has a blue matching of size  at least
 $|C_2|-\sqrt{\epsilon}t>(\frac{a}{2a+2}-\sqrt{\epsilon})t$.
Similarly, $H^b[C_3, C_4]$ has a blue matching of size at least $|C_4|-\sqrt{\epsilon}t>(\frac{a}{2a+2}-\sqrt{\epsilon})t$.
In total, we have that ${( {\cup_{i = 1}^4 {{C_i}} })^b}$ contains a
 blue matching $M_1$ with at least $(\frac{a}{a+1}-2\sqrt{\epsilon})t$ edges.

Note that all edges joining $C_5$ and $C_i$ are blue for $1\leq i\leq 4$, and each vertex  in $H$ has at most $\sqrt{\epsilon}t$ non-neighbors, so there is a blue fan with center $u\in C_5$ and blades in  $M_1$ of size at least
$$|E(M_1)|-\sqrt{\epsilon}t\geq \left(\frac{a}{a+1}-3\sqrt{\epsilon}\right)t
\overset{(\ref{eta-ep})}{>}\left(\frac{1}{2a+2}-0.05\beta\right)t$$
by noting $1/2\leq a<1$. This contradicts Claim \ref{large-fan-1}.

\medskip

 \noindent\textbf{Case B:} \ $\tau=4$

\medskip

Suppose that there exists a vertex $v_d\in C_1$ such that $\deg_{C^b_1}(v_d)\geq (\frac{a}{2a+2}-0.03\beta)t$, i.e., $v_d$ has at least $(\frac{a}{2a+2}-0.03\beta)t$ blue neighbors in $C_1$.
Let $C_1'$ be the  blue neighbors of $v_d$ in $C_1$.
We may assume that  $|C_1'|=(\frac{a}{2a+2}-0.03\beta)t$.
Recall that $|C_2|\geq \frac{a}{2a+2}t+1$ and all edges between $C_1'$ and $C_2$ are blue.
Thus for every $S\subseteq C_1'$,
\[
|N_{(H-K)^b}(S,C_2)|\geq |C_2|-\sqrt{\epsilon}t>\left(\frac{a}{2a+2}-\sqrt{\epsilon}\right)t> |S|.
\]
We apply Lemma \ref{hall} to the  subgraph induced by $C_1'\cup C_2$ in $H^b$
with $X=C_1'$ and  $Y=C_2$ to obtain a blue matching of size
$|C_1'|=(\frac{a}{2a+2}-0.03\beta)t$.
Since all (but at most $\epsilon t^2$) edges between $C_3$ and $C_4$  are blue, by a similar argument, we can find a blue matching of size
$(\frac{a}{2a+2}-\sqrt{\epsilon})t$   in $H^b[C_3,C_4]$.
Recall that every vertex in $H$ has at most $\sqrt{\epsilon}t$ non-neighbors
and  $1/2\leq a<1$, so we can find a matching of size at least
\begin{align*}
\left(\frac{a}{2a+2}-0.03\beta\right)t+\left(\frac{a}{2a+2}
-\sqrt{\epsilon}\right)t-\sqrt{\epsilon}t
\overset{(\ref{eta-ep})}{\geq}\left(\frac{a}{a+1}-0.05\beta\right)t\geq \left(\frac{1}{2a+2}-0.05\beta\right)t
\end{align*}
in $[N(v_d)]^b$, which together with $v_d$ forms a blue fan with blades more than $(\frac{1}{2a+2}-0.05\beta)t$.
This contradicts Claim \ref{large-fan-1}.

In the following, we may assume that
$\deg_{C^b_1}(v)< (\frac{a}{2a+2}-0.03\beta)t$ for every vertex $v\in C_1$.
We claim that $|C_1|\le(\frac{a}{a+1}-0.04\beta)t$. Otherwise, for every vertex $v\in C_1$,
\[\deg_{C^r_1}(v)\geq|C_1|-1-\sqrt{\epsilon}t-\deg_{C^b_1}(v)\geq|C_1|- \left(\frac{a}{2a+2}-0.02\beta\right)t>\frac{|C_1|}{2}.
\]
According to Lemma \ref{1-thm-1}, we obtain that $C^r_1$ hence $H^r$ contains a red cycle with more than $(\frac{a}{a+1}-0.3\beta)t$ vertices.
This contradicts Claim \ref{long-match-1}.
Note that $C_1$ is the largest red component, so we have
\begin{align}\label{ineq-11}
\frac{(1-\sqrt{\epsilon})t-1}{4}\leq|C_1|\leq\left(\frac{a}{a+1}-0.04\beta\right)t
\end{align}
 and
$| C_i|\geq \frac{a}{2a+2}t+1$ for every $2\leq i\leq 4$.
Note that
$|C_4|\leq (|C_1|+|C_2|+|C_3|+|C_4|)/4\leq t/4$, so we have
\[|C_1|+|C_4|\leq \left(\frac{a}{a+1}-0.04\beta\right)t+t/4=\left(\frac{5a+1}{4a+4}-0.04\beta\right)t.\]
Thus,
$|C_2|+|C_3|\geq |V(H)|-|K|-(|C_1|+|C_4|)\geq (\frac{3-a}{4a+4}+0.03\beta)t-1.$
It follows that
\begin{align}\label{ineq-10}
|C_2|> \left(\frac{3-a}{8a+8}+0.01\beta\right)t
\end{align}
 as $|C_2|\geq |C_3|$.
Since $|C_4|\leq t/4$, we have $|C_1|+|C_2|+|C_3|\ge 3t/4-\sqrt{\epsilon}t-1$. Therefore, we can take two disjoint subsets of $C_3$, say $C_3'$ and $C_3''$, such that \[|C_1|+|C_3'|=|C_2|+|C_3''|=
\left(\frac{1}{2a+2}-0.03\beta\right)t.
\]
By (\ref{ineq-11}) and (\ref{ineq-10}), we get that
\[
|C_3'|< \left(\frac{1-a}{4a+4}-0.03\beta+\frac{\sqrt{\epsilon}}{4}\right)t
\overset{(\ref{eta-ep})}{<}
\left(\frac{1-a}{4a+4}-0.02\beta\right)t
\]
and
\[
|C_3''|< \left(\frac{1}{2a+2}-0.03\beta\right)t-
\left(\frac{3-a}{8a+8}+0.01\beta\right)t
\overset{(\ref{eta-ep})}{<}\left(\frac{1}{8}-2\sqrt{\epsilon}\right)t.
\]
Therefore, $|C_3'|+\sqrt{\epsilon}t<|C_2|$ and $|C_3''|+\sqrt{\epsilon}t<|C_1|$.

%All edges between $C_4$ and $C_1\cup C_3'\cup C_2\cup C_3''$ are blue.
Denote $A=C_1\cup C_3'$ and $B=C_2\cup C_3''$.
Then the bipartite  graphs  $H^b[A,C_2]$ and $H^b[B,C_1]$ are almost blue complete bipartite graphs.

We claim that  the bipartite graph
$ H^b[A, B] $ contains a blue matching of cardinality at least $(\frac{1}{2a+2}-0.04\beta)t$.
Indeed, if $S\subseteq  C_3'\subseteq A$, then we have
\[|N_{(H-K)^b}(S, B)|\geq |C_2|-\sqrt{\epsilon}t>|C_3'|\geq |S|
\]
by noting that every vertex in $H$ has at most $\sqrt{\epsilon}t$ non-neighbors, and if $S\subseteq A$  and $S\cap C_1\neq {\O}$, then
\[
|N_{(H-K)^b}(S, B)|\geq |B|-\sqrt{\epsilon}t\geq|S|-\sqrt{\epsilon}t.
\]
Therefore, by Lemma \ref{hall}, the bipartite graph $H^b[A,B]$ and hence $H^b$
contains a blue matching  $M_2$ of cardinality at least $|A|-\sqrt{\epsilon}t=(\frac{1}{2a+2}-0.04\beta)t$ by (\ref{eta-ep}).
The claim follows.

Note that all (but at most $\epsilon t^2$) edges between $C_4$ and $A\cup B$ are blue, and every vertex in $H$ has at most $\sqrt{\epsilon}t$ non-neighbors.
Then there exists a vertex $w\in C_4$ whose blue neighborhood contains a
blue matching of cardinality  at least $|E(M_2)|-\sqrt{\epsilon}t=(\frac{1}{2a+2}-0.05\beta)t$ in $M_2$.
Thus we get a blue fan with center $w$ and at least $(\frac{1}{2a+2}-0.05\beta)t$ blades.
This leads to a contradiction by Claim \ref{large-fan-1}.

\medskip

\noindent\textbf{Case C:} \ $\tau=3$

\medskip
Note that all (but at most $\epsilon t^2$) edges joining $C_3$ and $C_i$ are blue for $1 \leq i \leq 2$. Since  each vertex in $H$ has at most $\sqrt{\epsilon}t$ non-neighbors and  $|C_1|\geq|C_2|\geq|C_3|\geq \frac{a}{2a+2}t+1$,
we may assume  that $|C_2|<(\frac{1}{2a+2}-0.02\beta)t$. Otherwise, by Lemma \ref{hall}, the bipartite graph
$H^b[C_1,C_2]$ contains a blue matching $M_3$ of size at least
$|C_2|-\sqrt{\epsilon}t \geq (\frac{1}{2a+2}-0.03\beta)t$ by noting $(\ref{eta-ep})$.
 Thus we can find a blue fan with center $v \in C_3$ and blades in $M_3$ of size at least
 $(\frac{1}{2a+2}-0.03\beta)t-\sqrt{\epsilon}t{\geq} (\frac{1}{2a+2}-0.04\beta)t$,
which contradicts Claim \ref{large-fan-1}. So $|C_2|<(\frac{1}{2a+2}-0.02\beta)t$ follows, and we have
\begin{align}\label{c-1}
\left(\frac{a}{a+1}+0.01\beta\right)t<t-\sqrt{\epsilon}t-1-2|C_2|\le|C_1|\leq \frac{t}{a+1}-2.
\end{align}
We claim that $C_1$ is 2-connected in $(H-K)^r$.  Otherwise, suppose that there exists
$K'\subseteq V(C_1)$ with $|K'|\leq 1$ such that $(C_1-K')^r$
is disconnected.  Since
$$\delta^r(C_1-K')\geq \delta^r(C_1)-1\overset{(\ref{ineq-9})}{\geq}
\frac{a}{2a+2}t-1,$$
we obtain that every red component in $(C_1-K')^r$ has size at least $\frac{a}{2a+2}t$,
thus $(C_1-K')^r$ has at most three red components by (\ref{c-1}) for every fixed $1/2\leq a<1$.
Thus $(H-K-K')^r$ has four  or five components, and so we are done from \textbf{Case A} or \textbf{Case B}.
Thus $C_1$ is 2-connected, which implies that the circumference
\[
c((C_1)^r)\geq \min\bigg\{2\delta^r(C_1), |C_1|\bigg\}\overset{(\ref{c-1}), (\ref{ineq-9})}{>}\frac{a}{a+1}t
\]
due to Lemma \ref{1-thm-3}.
This implies that  $C_1$ contains a red path with more than $(\frac{a}{a+1}-0.3\beta)t$ vertices in $H^r$, which leads to a contradiction by Claim \ref{long-match-1}.

\medskip

\noindent\textbf{Case D:}  \ $\tau=2$

\medskip

For this case, we claim that $C_1$ is  2-connected graph.
Otherwise, there exists a vertex set $K_1\subseteq C_1$ with $|K_1|\leq 1$
such that $(C_1-K_1)^r$ is disconnected.
Recall that $|C_1|\geq |C_2|\geq \frac{a}{2a+2}t+1$,
so we have
\[
\frac{(1-\sqrt{\epsilon})t-1}{2}\leq |C_1|\leq \frac{a+2}{2a+2}t-1.
\]
Moreover, it is clear that
$ \delta^r (C_1-K_1)\geq \delta^r (C_1)-1\geq \frac{a}{2a+2}t-1$ by (\ref{ineq-9}),
thus each  red component in  $(C_1-K_1)^r$ has size at least $\frac{a}{2a+2}t$.
This implies that $(C_1-K_1)^r$ contains at most four red components
for $1/2\leq a< 1$.
Therefore, $(H-K-K_1)^r$ contains either five, four or three red components, and so we are done from \textbf{Case A}, \textbf{Case B} or \textbf{Case C}.
Thus we conclude that  $C_1$ is  2-connected and $\delta^r(C_1)\geq \delta(H^r)-|K|\geq \frac{a}{2a+2}t$, which implies that $(C_1)^r$
 contains a red  path with more than $(\frac{a}{a+1}-0.3\beta)t$ vertices by Lemma \ref{1-thm-3}. This agian   contradicts Claim \ref{long-match-1}.

This completes the proof of Claim \ref{2-connect}.
\hfill$\Box$

\medskip
Now, by Lemma \ref{1-thm-3} and Claim \ref{2-connect}, we conclude that  the circumference
\[
c(H^r)\geq \min\{2\delta(H^r), |V(H^r)| \}
\geq \left(\frac{a}{a+1}t+2\right),
\]
where the last  inequality follows from Claim \ref{cla-match-1}.
Thus we obtain a red path with more than $(\frac{a}{a+1}t+2)$ vertices,
which contradicts Claim \ref{long-match-1}.

The proof of Part (I) is complete.
\hfill$\Box$

\bigskip

{\bf Part  (II)  \ $ a\geq 1$}

\bigskip
The lower bound $R(C_{2\lfloor an\rfloor}, F_n)\geq 4\lfloor an\rfloor-1$ is clear for every fixed $a\geq 1$.
Let $N=(4a+\gamma)n$, where $\gamma>0$ is a sufficiently small real number. Therefore, it suffices to show $R(C_{2\lfloor an\rfloor},F_n)\leq N$.
Thus we shall show that any red-blue
edge coloring of $K_N$ on vertex set $V$ yields either a red $C_{2\lfloor an\rfloor}$ or a blue $F_n$.
Suppose to the contrary that  for fixed $a\geq 1$ and large $n$, there exists a coloring that contains neither a red $C_{2\lfloor an\rfloor}$ nor a blue $F_n$. We aim to find a contradiction.

Similar as above, we apply the regularity lemma to obtain a partition of  $V$ with the corresponding properties, and $H$, $H^r$ and $H^b$ are defined similarly.
By a similar argument as Claim \ref{long-match-1} and Claim \ref{large-fan-1}, we get the following claims.
\begin{claim}\label{cla2-1}
$H^r$ contains no connected matching of size more than $\left(\frac{1}{4}-0.15\beta\right)t$.
\end{claim}

\begin{claim}\label{cla2-2}
$H^b$ contains no fan with at least $(\frac{1}{4a}-0.05\beta)t$ blades.
\end{claim}

We will also have the following claims.
\begin{claim}\label{cla2-3}
For each vertex $v\in V(H)$, $\deg_{H^r}(v)\geq \frac{2a-1}{4a}t+1$.
\end{claim}
\noindent{\bf Proof}.
On the contrary, we assume that $H^r$ contains a vertex $u$ such that
$\deg_{H^r}(u)\leq\frac{2a-1}{4a}t$. Since $\delta(H)\geq (1-2\sqrt{\epsilon})t-1$, we have
\[
\deg_{H^b}(u)\geq \left(1-2\sqrt{\epsilon}\right)t-\frac{2a-1}{4a}t-1
=\left(\frac{2a+1}{4a}-2\sqrt{\epsilon}\right)t-1.
\]
Denote $H_1=H[N_{H^b}(u)]$.
Note that every vertex in $H_1$ has at most $\sqrt{\epsilon}t$ non-neighbors.
Let $C_1$ and  $C_2$ be the vertex sets of  the largest blue  and red components in $H_1$ respectively.
Set $p=|C_1\setminus C_2|$.
By the same argument as Claim \ref{cla-match-1} step by step, we must have that $C_1=V(H_1)$ and
\begin{align*}
|C_1|=|V(H_1)|\geq \left(\frac{2a+1}{4a}-2\sqrt{\epsilon}\right)t-1, \;\; |C_2|>(1/2+2\sqrt{\epsilon})t, \;\; \text{and} \;\; p>20\sqrt{\epsilon}t.
\end{align*}
Then we have
\begin{align}\label{ineq-13}
|C_2|=|C_1|-|C_1\setminus C_2|\geq \left(\frac{2a+1}{4a}-2\sqrt{\epsilon}\right)t-p-1.
\end{align}

We first assume that $p=|C_1\setminus C_2|\geq (\frac{1}{4a}-0.05\beta)t$.
Note that $|C_2|>(1/2+2\sqrt{\epsilon})t$ and all (but at most $\epsilon t^2$) edges between  $C_1\setminus C_2$ and $C_2$ are blue.
Since each vertex in $H$ has at most $\sqrt{\epsilon}t$ non-neighbors,
we conclude that the bipartite graph $H^b_1[C_1\setminus C_2, C_2]$ contains a blue matching of size at least $(\frac{1}{4a}-0.05\beta)t$
by Lemma \ref{hall}.
Thus we can get a blue fan with center $u$ and at least
$(\frac{1}{4a}-0.05\beta)t$ blades, which contradicts Claim \ref{cla2-2}.

Thus we may assume $20\sqrt{\epsilon}t<p<(\frac{1}{4a}-0.05\beta)t$.
We  apply Lemma \ref{star-match} to the subgraph spanned by $C_2$ in $H_1$ with parameters $k=\sqrt{\epsilon}t+1$, $n_1=(\frac14-0.15\beta)t$, and $n_2=(\frac{1}{4a}-0.05\beta)t-p$ to obtain that
\[
R(S_k, n_1K_2, n_2K_2)=2n_1+n_2-1=\left(\frac{2a+1}{4a}-0.35\beta\right)t-p
-1\overset{(\ref{ineq-13}), (\ref{eta-ep})}{<}|C_2|.
\]
Since every vertex in $H_1$ has at most $\sqrt{\epsilon}t$ non-neighbors,
we can get a blue matching $M$ of size at least
$n_2=(\frac{1}{4a}-0.05\beta)t-p$  otherwise a red connected matching of size at least $(\frac14-0.15\beta)t$ will lead to a contradiction from Claim \ref{cla2-1}.
Note that
\begin{align*}
|C_2\setminus V(M)|= |C_1|-p-2n_2
\geq \left(\frac{2a-1}{4a}-2\sqrt{\epsilon}+0.1\beta\right)t+p-1
\overset{(\ref{eta-ep})}{>}
\left(\frac{2a-1}{4a}+6\sqrt{\epsilon}\right)t+p,
\end{align*}
so we have  $|C_2\setminus V(M)|>p+\sqrt{\epsilon}t=|C_1\setminus C_2|+\sqrt{\epsilon}t$.
Since all (but at most $\epsilon t^2$) edges between $C_1\setminus C_2$ and $C_2\setminus V(M)$ are blue,
according to the minimum degree of $H_1$, we obtain that
the bipartite graph $H^b_1[C_1\setminus C_2, C_2\setminus V(M)]$
contains a blue matching of size $p$ by Lemma \ref{hall}, which together with $M$ yield a blue
matching of size at least $n_2+p=(\frac{1}{4a}-0.05\beta)t$ in $H^b_1$.
Again, we can get  a blue fan with center $u$ and  at least
$(\frac{1}{4a}-0.05\beta)t$ blades in $H^b$, which contradicts Claim \ref{cla2-2}.
\hfill$\Box$

%By a similar analytic method as  Claim \ref{2-connect},  we can obtain the following claim. For simplicity, we omit the proof.

\begin{claim}\label{cla2-4}
$H^r$ is  2-connected.
\end{claim}

\noindent{\bf Proof}.
On the contrary,  there exists a subset $S\subset V(H)$ with $|S|\leq1$ such that $(H-S)^r$
is disconnected.
By Claim \ref{cla2-3}, we have
$$\delta^r(H-S)\geq \frac{2a-1}{4a}t,$$
and all red components of size at least $\frac{2a-1}{4a}t+1$.
Thus there are at most three red components in $(H-S)^r$.

If $(H-S)^r$ has three components $C_1,C_2$ and $C_3$ with $|C_1|\geq |C_2|\geq |C_3|$, then we have
\begin{align*}
|C_1|\geq |C_2|\geq |C_3|\geq \frac{2a-1}{4a}t+1,\;\; \text{and}\;\; \delta^r(C_i)\geq \frac{2a-1}{4a}t.
\end{align*}
By a similar argument as Case A of Part (I), $H^b$ contains a blue fan with more than $(\frac{1}{4a}-0.05\beta)t$ blades for $a\geq 1$.
This is a contradiction by Claim \ref{cla2-2}.

Therefore,  we  may assume that  $(H-S)^r$ has two components $C_1$ and $C_2$ with $|C_1|\geq |C_2|$.
It is clear that
\[
|C_2|\geq \frac{2a-1}{4a}t+1, \;\; \text{and} \;\;\frac{1-\sqrt{\epsilon}}{2}t-\frac12\leq|C_1|\leq \frac{2a+1}{4a}t-1.
\]
Suppose that there exists a vertex $u\in C_1$
such that $$\deg_{C^b_1}(u)\geq d:=\frac{|C_1|}{2}-\sqrt{\epsilon}t-1\geq \left(\frac{1}{4}-2\sqrt{\epsilon}\right)t.$$
Note that all (but at most $\epsilon t^2$) edges between $C_1$ and $C_2$ are blue and
each vertex in $H$ has at most $\sqrt{\epsilon}t$ non-neighbors.
Then $u$ together with $(\frac{1}{4}-2\sqrt{\epsilon})t$ blue neighbors in $C_1$
 and $(\frac{2a-1}{4a}-\sqrt{\epsilon})t$ blue neighbors in $C_2$ form a
 blue fan  with at least  $(\frac{1}{4a}-0.05\beta)t$ blades
by Lemma \ref{hall}.
This leads to a contradiction from Claim \ref{cla2-2}.
Thus we have $\deg_{C^b_1}(v)\leq d-1$ for every vertex $v\in C_1$.
It follows that
\[
\delta^r(C_1)\geq |C_1|-1-\sqrt{\epsilon}t-d+1
> \frac{|C_1|}{2}.
\]
By Lemma \ref{1-thm-1}, $(C_1)^r$ is pancyclic, which implies that $(C_1)^r$ contains a red cycle with length
\[
|C_1|\geq \frac{1-\sqrt{\epsilon}}{2}t-\frac{1}{2}
\overset{(\ref{eta-ep})}{>}\left(\frac{1}{2}-0.3\beta\right)t,
\]
which contradicts Claim \ref{cla2-1}.
\hfill$\Box$	

\medskip

Now note that  $(1-\sqrt{\epsilon})t\leq v(H)\leq t $ and $\delta(H^r)\geq \frac{2a-1}{4a}t+1$ from  Claim \ref{cla2-3},
it follows from Claim \ref{cla2-4} and Lemma \ref{1-thm-3} that the circumference
\[
c(H^r)\geq \min\{2\delta(H^r), |V(H^r)| \}
\geq \left(\frac{2a-1}{2a}t+2\right)>\left(\frac{1}{2}-0.3\beta\right)t,
\]
where the last inequality holds since $a\geq 1$.
This leads to a  contradiction from Claim \ref{cla2-1}.

The proof of Part (II) is complete.
\hfill$\Box$

\section{\textbf{Concluding remarks}}
In this paper, we are concerned with the asymptotic behavior of the Ramsey number $R(C_{2\lfloor an\rfloor}, F_n)$ when  $n$  is large and $a\geq 1/2 $ is fixed.
For fixed $0<a<1/2$ and large $n$, we also expect  to give a uniform  asymptotic behavior of $R(C_{2\lfloor an\rfloor}, F_n)$,
but  we encounter more obstacles for fixed $0<a<1/2$.
The graph  $G=3K_{2\lfloor an\rfloor-1}$  implies that $R(C_{2\lfloor an\rfloor}, F_n)\geq 6\lfloor an\rfloor-2$
for $2/5\leq a<1/2$. For $0<a<2/5$, the graph $G=K_{\lfloor an\rfloor-1}+\overline{K_{2n}}$
shows that $R(C_{2\lfloor an\rfloor}, F_n)\geq \lfloor an\rfloor +2n$.
When we consider the upper bound, our method will encounter more obstacles for
$0<a<1/2$. For example, the minimum  degree of $H^r$ maybe small and so we cannot conclude that $H^r$ is 2-connected for $0<a<1/2$.
Therefore, it would be interesting to determine the values of $R(C_{2\lfloor an\rfloor}, F_n)$ when $0<a<1/2$. Moreover, we do not know the behavior of $R(C_{2\lfloor an\rfloor+1}, F_n)$ for fixed $0<a\le 3/2$.

\bigskip

\noindent\textbf{Declaration of interests}

\medskip
The authors declare that they have no known competing financial interests or personal relationships that could have appeared to influence the work reported in this paper.

\end{spacing}

\end{document}